\magnification=1200
{\bf On the complexity of proper holomorphic mappings between balls}
\bigskip
John P. D'Angelo and Ji{\v r}\'\i\ Lebl,

Dept. of Mathematics,

University of Illinois,
 
1409 W. Green St.,

Urbana IL 61801

electronic addresses: jpda@math.uiuc.edu (D'Angelo) jlebl@math.uiuc.edu (Lebl)

\bigskip
{\bf Introduction}
\medskip

Let $B_n$ denote the unit ball in complex Euclidean space ${\bf C}^n$. Considerable literature is
devoted to the study of proper holomorphic mappings from $B_n$ to $B_N$. See for example [D1], [D2], [Fa1], [Fa2], [F1], [F2],
[H], [HJ], [HJX], [M], [P]. In this paper we make several
new contributions to this literature by establishing results linked to the notion of CR complexity theory.
Our results include a degree estimate for rational proper maps (Theorem 2), a new gap phenomenon for convex families of arbitrary
proper maps (Theorem 1), and an interesting result about inverse images (Theorem 4).

Two questions motivate this work:
how should one measure the complexity of a proper holomorphic mapping between $B_n$ to $B_N$, and how
are such measurements related to the domain and target dimensions? 
When a proper map $f$ between balls is smooth up to the boundary, and $n\ge 2$,
Forstneric [F1] proved that $f$ must be rational. For rational proper mappings, one 
possible approach to complexity is a degree estimate. For $n=1$, no such estimates are possible.
For $n=2$, a sharp bound is known in the monomial case; its proof involves an elaborate
graph-theoretic argument [DKR].  See [DLP] for results for $n\ge 3$. 
While a sharp bound for the degree of a rational proper mapping from $B_n$ to $B_N$ is not yet known, 
in Theorem 2 and Corollary 2 we give an improved (not sharp) bound.
To do so we combine a bound in the $2$-dimensional case proved by Meylan [M] with a technique developed by the authors and Peters in [DLP].
We obtain the following inequality: Assume $n\ge 2$. Let $f:B_n \to B_N$ be a rational proper mapping of degree $d$. Then
the degree $d$ of $f$ satisfies
$$   d \le {N(N-1) \over 2(2n-3)}. \eqno (0) $$

Although some of the work in this papers assumes rationality, it is natural
to weaken this assumption when possible. In Theorem 1 we prove a {\it gap phenomenon} for convex families of proper mappings 
with no assumptions on the regularity at the sphere of the map. 
We consider a one-dimensional convex family of (not necessarily rational) proper mappings from $B_n$ to $B_N$ that
preserves the origin. We prove, for $n\ge 2$, that $N$ must be at least $n+2$. 
Thus there is no one-dimensional convex family of origin preserving maps from $B_n$ to $B_{n+1}$ unless $n=1$.
By [Dor], there are many proper mappings in the codimension one case; our result shows that there is no convex family of them.
Convex families of proper mappings between balls arise naturally in this work; one reason relates to degree estimates.
Given a convex family of rational maps, only boundary members of this family are candidates for sharp degree estimates.
In fact the proof of the gap phenomenon is based upon passing to the endpoints of a one-dimensional family.

We define the crucial concept of convex family in Section II, but we give a short description now. A proper map $f$ from
ball to ball has a (vector-valued) power series that converges on compact subsets of the ball, and its squared Euclidean norm
defines a real-analytic function $z \to ||f(z)||^2$ there. The Taylor coefficients of $||f||^2$ at the origin 
determine a nonnegative definite Hermitian form.
The collection of such forms defines a convex cone and hence we take convex combinations of squared norms of proper maps.
Up to a linear map we can recover a holomorphic mapping from its squared norm. We therefore work with {\it norm equivalence};
proper maps $f$ and $g$ are norm equivalent if $||f||^2=||g||^2$ as functions. Convex families are thus the collection
of squared norms obtained by taking convex combinations. The simplest example (see [D], [D2], and [HJ])
is the {\it juxtaposition} of proper maps $f$ and $g$ from the same ball, but with possibly different targets. We let
$J_t(f,g)$ denote the norm-equivalence class of proper maps defined by 
$$ ||J_t(f,g)||^2 = {\rm cos}^2(t) ||f||^2 + {\rm sin}^2(t) ||g||^2. $$
This juxtaposition determines a one-dimensional convex family. The juxtaposition provides a homotopy
between an arbitary pair of proper maps from the same ball, as long as we allow sufficiently
many zero components and we identify $f$ with $f \oplus 0$ and $g$ with $g\oplus 0$. 

We next place known results into the context of complexity for proper mappings between balls.
Let $S^{2n-1}$ denote the unit sphere in ${\bf C}^n$ with its usual CR structure.
Assume $n\ge 2$ and let $f: S^{2n-1} \to S^{2N-1}$ be a CR mapping of class $C^{N-n+1}$.
By the result of Forstneric, $f$ must be the restriction of a rational mapping;
he also obtained a weaker bound than (0) for its degree.
As the codimension $N-n$ increases, the complexity of possible examples increases in a subtle fashion.
D'Angelo [D2] based on earlier joint work with Catlin [CD1] involving complex variables analogues
of Hilbert's 17th problem, established the following result.
Let $q$ be a polynomial on ${\bf C}^n$ that does not vanish on the closed unit ball.
Then there is an integer $N$ and a polynomial mapping $p: {\bf C}^n \to {\bf C}^N$ such that
${p \over q}$ is reduced to lowest terms and ${p \over q}$ maps $S^{2n-1} \to S^{2N-1}$.
The minimum possible target dimension $N$ is unbounded above, even when $n=2$ and the degree of $p$
is $2$. Thus every denominator is possible if we allow large enough target dimension. For $N < n$, however, the only maps are
are constant, and for $N=n$ the only nonconstant maps are linear fractional transformations.
For polynomial maps of unspecified degree, the dimension of the moduli space is arbitrarily large as $N$ tends to infinity.

\medskip 
{\bf Example 1}. Let $p=\sum c_{\alpha} z^\alpha$ be a polynomial of degree $d-1$ in $n$ variables with values in ${\bf C}^N$, 
and assume that the coefficient vectors all have sufficiently small Euclidean norm. In Proposition 3 we find a 
proper polynomial mapping $f$ of degree $d$ whose jet of order $d-1$ is $p$. The mapping $f$
is determined up to norm equivalence by the collection of inner products $\langle c_\alpha, c_\beta \rangle$. 
The possible values of these inner products can be regarded as parameters. The dimension of the parameter space over the real numbers,
determined precisely in Proposition 3, obviously tends to infinity as $d$ does. 

\medskip

We will also give some restrictions on the dimensions of convex families of rational proper mappings 
from $B_n$ to some $B_N$ in terms of $n,N$ and the degree $d$.
Theorem 3 applies only in the monomial case. We show that degree bounds for monomial proper maps
can be extended to degree bounds for monomial families. In particular, if we have a bound $d\le c(n,N)$ for the degree $d$ in terms
of the domain and target dimensions, then we obtain the bound $d\le c(n,N-k)$ for the degree of each element in a $k$-dimensional family.
We show by example that the result fails in the polynomial case.

The inverse image of a point under a proper holomorphic mapping is both compact and a complex analytic variety, and hence it is a finite set.
In Theorem 4 we provide additional new information. For example we show that
if $p$ is a proper polynomial mapping between balls, then as a set
either $p^{-1}(0)$ is empty or $p^{-1}(0) = \{0\}$. We give two proofs.
In the rational case the inverse image can be any finite set of points in the domain ball, as long as 
we allow the degree of the map to be sufficiently large, in which case the target ball must be
of sufficiently high dimension. To get control on the situation, we fix the denominator $q$ of a rational $f$. 
From $q$ we determine a finite set of candidates for the elements of $f^{-1}(0)$.
The proof combines polarization and homogenization techniques. 
In the polynomial case, while $p$ need not be injective, 
the origin itself cannot have multiple inverse images; the various inverse images of nearby points must coalesce there. 
Theorem 4 seems to be a first step toward a deeper algebraic understanding of the rational case.

In the final section we provide subtle examples of convex families
and proofs of various restrictions on the four integers $n,N,d,k$
describing the complexity: the domain dimension, the target dimension, the degree, and the dimension of the family.

The first author acknowledges support from NSF grant DMS 05-00765. 
The authors began discussing this kind of problem at MSRI in 2005;
the first author ran a graduate student
experience in CR Geometry which the second author attended while a student.
Both acknowledge support from MSRI. Both also acknowledge AIM for the meeting on Complexity Theory
in CR Geometry held in 2007. Both authors also thank Han Peters, Xiaojun Huang, and Shanyu Ji
for helpful discussions.

\medskip
{\bf I. Preliminaries}
\medskip

The condition for $f$ being a proper mapping from $B_n$ to $B_N$ 
is of course that $||f(z)||^2$ tends to unity whenever $||z||^2$ tends to unity from below.
We will therefore often work with $||f||^2$ rather than with $f$ itself.
Let $\psi$ denote an inclusion map of a ball into the equator of a larger dimensional ball. If $g$ is proper between balls, then
so is $\psi \circ g$; of course $||g||^2 = ||\psi \circ g||^2$. We therefore avoid some nuisances by considering 
two mappings from the same ball to be equivalent if their squared norms are the same function. 
Hence we will say that proper maps $f$ and $g$ from the same ball to possibly different balls
are {\it norm-equivalent} if $||f||^2=||g||^2$ as functions.

Observe that the automorphism group of the ball is transitive. Hence we may assume after composition with an automorphism
that a proper map satisfies $f(0)=0$. We make this assumption, unless stated otherwise, whenever $f$ is a proper mapping between balls.

Let $z \to R(z,{\overline z})$ be a real-analytic real-valued
function defined on $B_n$; we may expand it in a power series:
$$ R(z,{\overline z})  = \sum_{\alpha,\beta} c_{\alpha \beta} z^\alpha {\overline z}^\beta.$$
We call $(c_{\alpha \beta})$ the {\it underlying Hermitian form} corresponding to $R$.
Following [DV] we say that $R \in {\cal P}_\infty$ if there is a holomorphic function $f$ with values in a Hilbert space 
such that $R(z, {\overline z}) = ||f(z)||^2$. In this paper the underlying Hermitian form will always have finite rank, and hence we may assume that the Hilbert space is finite-dimensional. The simplest test for $R \in {\cal P}_\infty$ is that 
the underlying Hermitian form be nonnegative definite.

We next make the connection to proper mappings between balls.
Given such a mapping $f$, we may expand it and its squared norm in Taylor series about the origin. Thus for elements
$c_\alpha \in {\bf C}^N$ we have:
$$ f(z) = \sum_\alpha c_\alpha z^\alpha \eqno (1) $$
$$ ||f(z)||^2 = \sum_{\alpha,\beta} \langle c_\alpha, c_\beta \rangle z^\alpha {\overline z}^\beta. \eqno (2) $$
Thus $c_{\alpha \beta} = \langle c_\alpha, c_\beta \rangle$ and hence this form is nonnegative definite.
We define the {\it rank} ${\cal R}(f)$ to be the rank of its underlying Hermitian form. Then
${\cal R}(f)$ is the smallest integer $N_0$ for which there is a map $g$,
norm-equivalent to $f$ and mapping $B_n$ properly to $B_{N_0}$.

\medskip
{\bf Remark 1}.  Let $f : B_n \to B_N$ be a proper mapping with $f(0)=0$. Then the embedding dimension 
of $f$ equals ${\cal R}(f)$. 

Proof: By linear algebra, in all cases the embedding dimension of $f$ is at most ${\cal R}(f)$. 
To prove the reverse inequality we must assume that $f(0)=0$.
Suppose the embedding dimension of $f$ is $N-k$.
Then there exists an automorphism $\phi$ such that $k$ components of $\phi \circ f$ vanish.
The numerator of these components defines an affine function of rank $k$ applied to $f$. Because
$f(0)=0$, this affine function is linear. Therefore the rank of $f$ is at most $N-k$.
\ $\spadesuit$

\medskip

Observe that a convex combination of nonnegative definite Hermitian forms is also
nonnegative definite. More generally we have the following result:

\medskip
{\bf Lemma 1}. For $1 \le j\le k+1$, let $f_j$  be a proper mapping from $B_n$ to $B_{N_j}$. For $\lambda \in {\bf R}^k$, define $R_\lambda$ by
$$ R_\lambda = \sum_{j=1}^k \lambda_j ||f_j||^2 + (1- \sum _{j=1}^k \lambda_j) ||f_{k+1}||^2. \eqno (3) $$ 
Let $K$ denote the set of $\lambda \in {\bf R}^k$ for which $R_\lambda \in {\cal P}_\infty$. Then $K$ is compact and convex.

\medskip

For $k\ge 2$, the set $K$ from Lemma 1 can be complicated. See Section V.
Even when $k=1$, it can strictly contain $[0,1]$.  We give a simple example.

\medskip
{\bf Example 2}. Put $||f||^2 = a ||z||^4 + (1-a) ||z||^2$ and put $||g||^2 = b||z||^4 + (1-b)||z||^2$, and assume that $0\le b < a \le 1$. 
Let $R_\lambda$ be a convex combination:
$$ R_\lambda = \lambda ||f||^2 + (1-\lambda)||g||^2 = (\lambda a + (1-\lambda)b) ||z||^4 + (\lambda(1-a) +(1-\lambda)(1-b))||z||^2. \eqno (4) $$
It follows from (4) that $R_\lambda \in {\cal P}_\infty$ if and only if 
$$ {-b \over a-b} \le \lambda \le {1-b \over a-b}. \eqno (5) $$
Thus the set of $\lambda$ for which $R_\lambda \in {\cal P}_\infty$ is a closed bounded interval strictly containing $[0,1]$.
\medskip

\medskip
{\bf II. Families}
\medskip
We next discuss several possible notions of {\it family} of proper maps. The simplest notion, but one too general for our purposes, would be to consider continuous maps from the parameter space to a set indexed by the collection of multi-indices in $n$ variables. Thus we would
demand, for each $\alpha$, that $\lambda \to c_\alpha(\lambda)$ be continuous, and that
$$ f(z) = \sum_{\alpha} c_\alpha(\lambda)  z^\alpha $$
be a proper holomorphic map from $B_n$ to $B_N$.

We give a simple example. Let the parameter space be the unit disk. For $n=N=1$ we define
$c_0(\zeta) = -\zeta$, and for $a \ge 1$, we define $c_a(\zeta) = {\overline \zeta}^{a-1} (1 - |\zeta|^2)$. Then each map
$\zeta \to c_a(\zeta)$ is continuous, and for $|\zeta| < 1$, the map $\phi_\zeta$ defined by
$$ \phi_\zeta(z) = \sum c_a(\zeta) z^a = {z-\zeta \over 1 - {\overline \zeta}z} $$
defines a proper holomorphic mapping (in fact an automorphism) from $B_1$ to itself. This particular family is not closed under limits;
as $|\zeta|$ tends to $1$, the limit map is a constant, and hence not proper.

Consider in general a proper map $f$ between balls, and let $u$ and $v$ be
automomorphisms of the domain and target balls. Then $g=v \circ f \circ u$ is also a proper map, and we say that $f$ and $g$ are {\it spherically equivalent}. An automorphism is a linear fractional transformation
determined by two pieces of data: the image of the origin and a unitary mapping. Spherical equivalence differs from norm equivalence, as Example 3 shows.
One can also combine the notions. Assume $f$ and $g$ are maps from the same ball but with different target balls.
It is possible for $f$ to be spherically equivalent,
via nonlinear automorphisms, to $g \oplus 0$. Then $f$ and $g$ are neither norm equivalent nor spherically equivalent, yet they define maps
whose properties are the same for many purposes.

\medskip
{\bf Example 3}. Put $n=2$. Define $f$ and $g$ by $f(z)= (z_1,z_1z_2,z_2^2)$
and $g(z) = (z_1^2, z_1z_2,z_2)$. Then $f=g\circ u$, where $u$ interchanges the variables. Hence $f$ and $g$ are spherically
equivalent, but they are not norm equivalent. On the other hand, the maps $z$ and $z \oplus 0$ are norm equivalent but not spherically
equivalent because they have different target dimensions.
\medskip

The clue for our definition of convex families comes from [D2]. Consider proper mappings $f$ and $g$ with the same domain ball, 
but with targets $B_{N_1}$ and $B_{N_2}$. We may regard them as homotopic in the following fashion. 
For $t \in [0, {\pi \over 2}]$, define $J_t(f,g)$ by 
$$ J_t(f,g) = {\rm cos}(t)f \oplus {\rm sin}(t) g.  \eqno (6) $$
Then $ J_t(f,g)$ defines a proper map to $B_{N_1 + N_2}$.
If we identify $f$ with $f \oplus 0$ and $g$ with $0 \oplus g$, then the family $J_t(f,g)$ defines a homotopy between
$f$ and $g$. (The homotopy parameter is ${\rm sin}(t)$ rather than $t$ itself.)
Of course $f$ and $f \oplus 0$ are norm-equivalent. Furthermore $J_t(f,g)$ is norm equivalent to $J_t(g,f)$.
Computing the squared norm of $J_t(f,g)$ yields
$$ ||J_t(f,g)||^2 = {\rm cos}^2(t)||f||^2 + {\rm sin}^2(t) ||g||^2. $$

 Replacing ${\rm cos}^2(t)$ by $\lambda$ shows that
$||J_t(f,g)||^2$ is a convex combination of the squared norms $||f||^2$ and $||g||^2$ and that the Hermitian matrix
$c_{\alpha \beta}(t)$ corresponding to $J_t(f,g)$ depends linearly on $\lambda$. We
obtain a $k$-{\it dimensional convex family} by iterating this operation $k$ times. We make the following slightly more general
definition.

\medskip
{\bf Definition 1}. A zero-dimensional family of proper mappings from $B_n$ is a single proper map $f:B_n \to B_N$ for some $N$.
Let $f_i: B_n \to B_{N_i}$ be proper mappings whose squared norms define linearly independent Hermitian forms. For $\lambda \in {\bf R}^k$, define a collection of real analytic functions $R_\lambda$ by
$$ R_\lambda = \sum_{j=1}^k \lambda_j ||f_j||^2 + (1- \sum_j \lambda_j) ||f_{k+1}||^2. \eqno (7) $$
Let $K$ be the largest convex set in ${\bf R}^k$ for which $\lambda \in K$ implies $R_\lambda \in {\cal P}_\infty$.
The $k$-{\it dimensional convex family} ${\cal F}$ consists of
the proper mappings $f$ for which
$ ||f||^2 = R_\lambda$ for some $\lambda \in K$. The {\it rank} of a family is the maximum of the ranks of its members.

\medskip
The expression in (7) can be regarded as a linear combination of nonnegative definite
Hermitian forms, where the sum of the coefficients equals $1$. Note therefore that the dimension of the family
is $k$ while the Hermitian forms corresponding to the squared norms span a $k+1$-dimensional subspace.
It is obvious that a convex combination, where the coefficients are nonnegative, of such forms is also nonnegative definite. Here we allow some of the coefficients to be negative, as long as the form defined by (7) is nonnegative definite. Example 2 provided a case where
a coefficient could be negative. In general the set $K$ from Lemma 1 can be rather complicated. See the discussion after Proposition 3.

When ${\cal F}$ is a family we will write $f_\lambda$ for its members. Thus
$ ||f_\lambda||^2 = R_\lambda,$ where $f_\lambda$ is determined up to norm equivalence. The index $\lambda$ does not mean a
component of $f$.

\medskip {\bf Example 4}. Let $f,g,h$ be proper mappings from the same ball. 
Define $J_t(f,g,h)$ by 
$$ {\rm cos}(t_1) {\rm cos}(t_2) f \oplus {\rm sin}(t_1) {\rm cos}(t_2)g  \oplus {\rm sin}(t_2) h. \eqno (8) $$
Computing squared norms in (8) we obtain
$$ ||J_t(f,g,h)||^2 = {\rm cos}^2(t_1) {\rm cos}^2(t_2)||f||^2  + {\rm sin}^2 (t_1) {\rm cos}^2 (t_2)||g||^2  + {\rm sin}^2 (t_2) ||h||^2. \eqno (9) $$
Thus $\lambda_1 = {\rm cos}^2(t_1){\rm cos}^2(t_2)$ and $\lambda_2 = {\rm sin}^2 (t_1) {\rm cos}^2 (t_2)$. 
We began by linking $f$ and $g$ as members of a $1$-dimensional family, and then we linked $h$ to these maps.

\medskip
To make Example 4 more specific, we take for example $h(z)=z$, $g(z) = z \otimes z$, and $f(z)= z\otimes z \otimes z$.
We obtain a two-dimensional convex family of cubic polynomial mappings. We next provide further examples of convex families
and illustrate one difference between the one-dimensional situation and the higher dimensional case.
\medskip
{\bf Example 5}. The mapping $f$ defined by 
$$f(z)= (z_1,{\sqrt {2} \over 2}z_1z_2, {\sqrt {2} \over 2} z_2^2, {\sqrt {2} \over 2} z_2)$$ 
is proper from $B_2$ to $B_4$. Note that $f$ is in the interior of the $1$-dimensional family 
$$ (z_1, {\rm cos}(t)z_1z_2, {\rm cos}(t)z_2^2 ,{\rm sin}(t)z_2^2). \eqno (10)$$ 

We prove a surprising result about one-dimensional convex families. Even if we allow
proper mappings with no regularity assumption, there are restrictions on the target dimension for families. We begin
with a simple example which indicates why the case of domain dimension $1$ is special.
The crucial point, when $n=N=1$, is that proper mappings need not be automorphisms. For $n=N\ge 2$,
proper self maps of $B_n$ are necessarily automorphisms. [P]

\medskip
{\bf Example 6}. For $n=1$ put $f_\lambda(z) = (a z, bz^2)$ where $|a|^2 =\lambda$ and $|b|^2 =1-\lambda$. Then $f_\lambda$
defines a one-dimensional convex family of proper mappings from $B_1$ to $B_2$, and $f_\lambda(0)=0$.
For $0< \lambda <1$ the map lies in the interior of that family.

\medskip
The next result shows for higher domain dimensions that no one-dimensional convex family exists in codimension one.
By [Dor], there are many proper mappings in the codimension one case; we show that no convex families exist.

\medskip
{\bf Theorem 1}. Let ${\cal F}$ be a positive-dimensional convex family of proper maps
from $B_n$ to $B_N$. Assume $n\ge 2$ and $f(0) = 0$ for all $f \in {\cal F}$.
Then $N\ge n+2$. Thus for $n\ge 2$ there is no one-dimensional convex family of origin-preserving
proper maps from $B_n$ to $B_{n+1}$.

Proof. Let ${\cal F}$ be a convex family of maps from $B_n$ to $B_N$. Let $c_{\alpha \beta}(\lambda)$ be
the underlying Hermitian form corresponding to each $f_\lambda$. According to Lemma 1, we may extend the family by allowing $\lambda$
to live in the maximal compact convex set $K$ for which $R_\lambda$ is a squared norm for each $\lambda$. Consider any boundary
point $\mu$ of $K$. When an eigenvalue of $c_{\alpha \beta}(\mu)$ is positive, the same eignvalue must be 
positive in a neighborhood of $\mu$ because the eigenvalues are continuous in $\lambda$.
Since $\mu$ is on the boundary of $K$, it must be true that at least one eigenvalue of $c_{\alpha \beta}(\mu)$ which was positive
in the interior of $K$ will vanish at $\mu$.
If in addition the parameter space has dimension one, then $K$ is an interval $[a,b]$.  Thus an eigenvalue 
will vanish at each endpoint of that interval. If $N=n+1$, then each of the endpoint maps will have rank at most
$n$. If the rank were less than $n$, then the map would be a constant. Since $f_\lambda(0)=0$ for all $\lambda$, this case doesn't arise.
Therefore each of the endpoint maps must have rank $n$, but the only proper holomorphic self mappings of a ball (for $n\ge 2$)
are automorphisms. Since the origin is preserved, the endpoint maps must each be unitary linear transformations and hence
have squared norm $||z||^2$. It follows that any convex combination of them has squared norm $||z||^2$, and hence the family is actually
zero-dimensional. Therefore the case $N=n+1$ is ruled out when $n\ge 2$. We conclude that $N\ge n+2$. \ $\spadesuit$ 

\medskip 
{\bf Corollary 1}. Assume $n \ge 2$. Let $f : B_n \to B_N$ and $g : B_n \to B_M$ be proper maps each preserving the origin. 
Suppose for some $t \in (0,1)$ that the juxtaposition $J_t(f,g)$ is spherically equivalent to a map of the form $h\oplus 0$, where
$h : B_n \to B_{n+1}$. Then both $f$ and $g$ are linear. 

Proof. Consider the maximal one-dimensional family ${\cal F}$
containing $J_t(f,g)$. By the hypothesis and Remark 1, there is an interior element of ${\cal F}$ with rank at most $n+1$, 
and hence the rank of all elements in ${\cal F}$ is at most $n+1$. By Theorem 1, the rank must 
actually be $n$, in which case $f$ and $g$ are linear and the family is zero-dimensional.
\  $\spadesuit$

\medskip
If $n\ge 2$ and we do not assume $f(0)=0$, then $N=n+1$ is possible in Theorem 1.
If ${\cal F}$ is a $1$-dimensional family mapping to the $n+1$ ball, then we must have
$ f_\lambda = {\sqrt \lambda} u \oplus \sqrt{1-\lambda^2}$ for some automorphism $u$ of $B_n$.
When we assume $f(0)=0$, we must have $\lambda=1$ and the automorphism must be unitary. When $n=1$, Example 6 shows that
a $1$-dimensional family preserving the origin is possible.

\medskip {\bf Remark 2}. Theorem 1 is sharp. Example 5 gives a $1$-dimensional family of quadratic polynomial
mappings from $B_2$ to $B_4$. See [D] and [HJ] for uses of this particular family.

\medskip

{\bf III. Degree estimates for rational maps}
\medskip

The first result in this section uses a pullback procedure introduced in [DLP]. 
We then prove a related result for monomial families. We close by showing in the rational case that
interior elements in a convex family cannot give sharp degree estimates.

\medskip
{\bf Theorem 2}. Suppose $N\ge 2$, and the estimate $d\le c(2,N)$ holds for the degree of rational proper mappings from $B_2$
to $B_N$. Then the estimate 
$$d \le {c(2,N) \over 2n-3}$$
 holds for the degree of rational mappings from $B_n$ to $B_N$.

Proof. For $n\ge 2$ consider the map $V:B_2 \to B_n$ defined by 
$$ V(z,w) = (z^{2n-3}, c_1 z^{2n-5}w, ..., c_s z^{2n-3-2s} w^s, ... , c_{n-2} z w^{n-2}, w^{2n-3}). \eqno (11) $$
With the correct choice of constants (See [D] or [DLP]) the map $V$ maps $B_2$ properly to $B_n$; the result is a monomial
proper mapping invariant under a certain representation of a cyclic group of order $2n-3$. Given any rational
proper mapping $g$ of degree $d$ from $B_n$ to $B_N$, we may perform a unitary transformation to ensure that
the first component includes a monomial of degree $d$. In the pullback that monomial gets raised to the power
$2n-3$. Since all the terms in $V$ are of degree at most $2n-3$ we can ensure that the degree of $g \circ V$ is 
in fact $(2n-3)d$. It follows that $(2n-3)d = {\rm deg} (g \circ V)$. 

The composition $g \circ V$ maps $B_2$ to $B_N$ and we therefore have
$$ (2n-3)d = {\rm deg} (g \circ V) \le c(2,N). $$
Hence $d \le { c(2,N) \over 2n-3}$. $ \ \spadesuit$

\medskip

Although Theorem 2 is not sharp, it is useful. For example, Meylan [M] has established the bound
$d \le {N(N-1) \over 2}$. Combining her bound with Theorem 1 yields the following:

\medskip
{\bf Corollary 2}. Let $f:B_n \to B_N$ be a rational proper mapping of degree $d$ and assume $N\ge 2$. Then
the degree $d$ of $f$ satisfies
$$   d \le {N(N-1) \over 2(2n-3)}. $$
Proof. By the above results, 
$$ d \le { c(2,N) \over 2n-3} \le {N(N-1) \over 2(2n-3)}. \ \spadesuit $$

Our next result enables us to extend bounds for the degree of a monomial
proper mapping to bounds for the degree of proper mappings in a $k$-dimensional convex family.
It does not hold for polynomial families in general.

Let $f$ be a monomial proper mapping between balls. The squared norm of $f$ depends on only the variables $|z_i|^2$,
and we express things more simply by writing 
$$ x=(x_1,...,x_n) = (|z_1|^2,...,|z_n|^2). $$
As in [DLP] we write 
$$ p(x) = p(x_1,...,x_n)= p(|z_1|^2,...,|z_n|^2)= ||f(z)||^2, \eqno (12) $$
where $p$ is a polynomial in $x$ whose coefficients are nonnegative.
We call $p$ the {\it real form} of $f$. The condition for being proper is simply that $p(x)=1$ on the hyperplane
given by $\sum x_j=1$. The number of terms in $p$ equals the rank ${\cal R}(f)$. 

Consider any bound $d \le c(n,N)$ for the degree of a proper rational
map $f:B_n \to B_N$. Since the map to $B_{N+1}$ defined by $z \to (f(z), 0)$ 
also has degree $d$, we may assume that $N \to c(n,N)$ is a nondecreasing function.

\medskip
{\bf Theorem 3}. Let ${\cal F}$ be a convex family of dimension $k$ consisting
of monomial proper maps from $B_n$ to $B_N$, where $n \ge 2$.
Assume the bound $d\le c(n,N)$ holds for all monomial maps of degree $d$ from $B_n$ to $B_N$. Then, for all
$f \in {\cal F}$, we have ${\rm deg}(f) \le c(n,N-k)$.

Proof. We proceed by induction, with the basis step $k=0$ being the hypothesis. Suppose the result holds
for $k-1$ dimensional families. Thus we assume for every $k-1$-dimensional family ${\cal F}$ of monomial maps
that 
$${\rm deg}(f) \le c(n,N-(k-1))$$
for every $f \in {\cal F}$. We may assume that $N \to c(n,N)$ is nondecreasing.
Suppose $g_\lambda$ lies in a $k$-dimensional family.
Let $p_\lambda(x)$ be the real form of $g_\lambda$ as in (12). We obtain
$$ p_\lambda(x) = \sum c_\alpha(\lambda)x^\alpha, $$
where $\lambda \to c_\alpha(\lambda)$ is an affine function for each multi-index $\alpha$.
As in Definition 1, let $K$ be the largest convex set
such that $\lambda \in K$ implies that $c_\alpha(\lambda) \ge 0$ for all $\alpha$. On the boundary of $K$ we must have 
$c_\alpha(\lambda)=0$ for at least one $\alpha$. We solve this affine equation to eliminate some $\lambda_j$, and substitute
this value in the formula for $p_\lambda$ to obtain $p_{\lambda_0}$. Thus 
$$ N= N(p_\lambda) \ge N(p_{\lambda_0})+1 = N_0 +1, \eqno (13) $$
where $N(p)$ denotes the number of terms in $p$. The collection
of these maps $p_{\lambda_0}$ defines a $k-1$ dimensional convex family of monomial mappings. Obviously 
${\rm deg}(p_{\lambda_0}) \le {\rm deg}(p_\lambda)= d$, but we need equality. We
claim that we can always choose the index $\alpha$ to guarantee that ${\rm deg}(p_{\lambda_0}) = {\rm deg}(p_\lambda)$. If this claim were false, then every boundary point of $K$ would be of degree less than $d$.  Since $K$ is convex, no interior point
could be of degree $d$. Therefore we may assume that ${\rm deg}(p_{\lambda_0}) ={\rm deg}(p_\lambda)$.

Now using the induction hypothesis and (13) we obtain
$$ {\rm deg}(p_\lambda) = {\rm deg}(p_{\lambda_0}) \le c(n,N_0-(k-1)) \le c(n,N-k), $$
providing the induction step. In this argument we have eliminated one of the $\lambda_j$, maintained the degree, and 
lowered the dimension of the family by one. \ $\spadesuit$

\medskip

The sharp bound $d \le 2N-3$ holds for monomial proper mappings from $B_2$ to $B_N$. The bound
$d \le {4 \over 3} {2N-3 \over 2n-3}$ holds for general $n$ but it is not sharp. For $n$ large compared with $d$ and at least $3$,
the bound $d \le {N-1 \over n-1}$ holds. In fact, examples show that the floor of ${N-1 \over n-1}$ gives the sharp bound.
See [DKR] and [DLP]. Combining these bounds with Theorem 3 yields
bounds for the dimensions of convex families. We state one such result.

\medskip
{\bf Corollary 3}. Let ${\cal F}$ be a $k$-dimensional convex family of monomial maps from $B_n$ to $B_N$. Then, for each $f \in {\cal F}$,
$$ {\rm deg}(f) \le {4 \over 3} {2(N-k)-3 \over 2n-3}. \eqno (14) $$

In Theorem 3 we were able to replace $N$ by $N-k$ in the bound $c(n,N)$. We cannot do so for polynomial or rational maps except in
the case $k=1$. See Proposition 1 below. Thus Theorem 3 is restricted to monomial mappings.
We provide an example where Theorem 3 fails in the polynomial case.

\medskip {\bf Example 7}. Let $n=2$ and write the variables as $(z,w)$.
Consider the norm equivalence class of a polynomial mapping $f:B_2 \to B_5$ of the form
$$ f(z,w) = Az + B w + C z^2 + Dzw + Ew^2. \eqno (15) $$
Each coefficient in (15) is an element of ${\bf C}^5$; by an elementary calculation (See page 168 of [D] for example), 
we may specify the linear coefficients arbitrarily as long as they are sufficiently small, and then the quadratic
coefficients are determined. The condition {\it sufficiently small} can be stated $I - L^*L$ is positive semi-definite,
where $L$ is the linear map defined by $L(z,w)= Az+Bw$.
The relevant parameters become the inner products $||A||^2$, $||B||^2$, and $\langle A, B\rangle$. 
These numbers give $4$ real parameters, as the squared norms are nonnegative. The underlying Hermitian matrix as in (2)
is of size five by five; its top two-by-two block is formed from these three parameters, and the remaining entries
are determined by them. The space of such two-by-two blocks is also $4$ real dimensional. It follows
from Definition 1 that $||f||^2$ determines a $4$-dimensional family of polynomial mappings of degree $2$ from $B_2$ to $B_5$. 
For any possible sharp bound $c(n,N)$ we must have $c(2,1)=0$, and hence
the conclusion $2= d \le c(n,N-k) = c(2,1)=0$ fails.
\medskip
In the proofs of Theorems 1 and 3 we passed to the boundary of the family. The same technique yields the following result, which
generalizes a result from [L].

\medskip {\bf Proposition 1}. Let ${\cal F}$ be a convex family of rational
proper holomorphic mappings from $B_n$ such that
${\cal F}$ is positive dimensional, has rank $r$, and is of generic degree $d$.
Let $c(n,N)$ denote an upper bound for the degree of a rational proper mapping from $B_n$ to $B_N$.
Then $d \le c(n,r-1)$.

Proof. Let $K$ be the compact convex set from Lemma 1. By the convexity of $K$, there must be some $f$ of degree $d$
and in the boundary of the family ${\cal F}$. The condition ``being of rank $r$'' is an open condition. By therefore
passing to an appropriate boundary point of the family, 
we obtain a rational mapping $f$ of the same degree whose rank is at most $r-1$. 
Thus we may assume that $f:B_n \to B_{r-1}$, and hence $d\le c(n,r-1)$. \ $\spadesuit$

\medskip
{\bf IV. Inverse images of the origin}
\medskip

The main result in this section considers the possibilities for the inverse
image of $0$ under a rational proper mapping $f$ between balls with given denominator $q$.
We construct a finite set $S(q)$ which contains $f^{-1}(0)$.
In the polynomial case this finite set is the origin alone,
and hence we conclude that either $f^{-1}(0)$ is empty or it consists of $0$ alone.  
We have been following the convention, without loss of generality, that $f(0)=0$ when $f$ is a 
proper mapping between balls. For polynomials,
Theorem 4 reveals that $0$ plays a special role. Since a proper map is a {\it finite} map,
the inverse image of $0$ must be a finite set, together with multiplicities. There is no other restriction
on this inverse image in the rational case, as the following standard result [D] shows.

\medskip
{\bf Proposition 2}. Let $a_1,...,a_k$ be an arbitary finite collection of points in $B_n$ (repeats are allowed).
Then there is an $N$ and a proper rational mapping $f:B_n \to B_N$ such that
$f^{-1}(0)= \{ a_1,...,a_k \}$.

Proof. Let $\phi_{a_j}$ be an automorphism of $B_n$ with $\phi_{a_j} (a_j) = 0$. Let $f$ denote the tensor product
of these automorphisms. Then $f$ is a proper rational mapping from $B_n$ to $B_N$, where $N=n^k$. Furthermore
its zero set consists precisely of the points $a_j$. \ $\spadesuit$

\medskip
Let $q$ be a polynomial that does not vanish on the closed unit ball. By [D2], $q$ is the denominator of 
some proper rational mapping $f$ between balls that is reduced to lowest terms. It is easy to see that $f^{-1}(0)$ 
is a subset of the reflection across the unit sphere of the variety defined by $q$. We will narrow down
this possibility to a finite set. Let us assume that $q$ is of degree $d-1$.

Write $q = \sum _{j=0}^{d-1} q_j$ as 
its expansion in terms of homogeneous parts. Using an idea from [D4] we define a homogeneous polynomial $Hq$ of degree $d$ by

$$ Hq(w, z) = \sum_{j=0}^{d-1} \langle w, z \rangle^{d-j} q_j(w).  \eqno (16) $$
Note that the sum is divisible by $\langle w,z \rangle$.
Then $Hq(w,z)$ is homogeneous in $w$ of degree $d$, and it is an anti-holomorphic polynomial in $z$.
We claim that the set $S(q)$ of $z$ for which $H(w,z)$ vanishes identically is finite, from which we
we obtain Theorem 4 on inverse images. First we prove that $S(q)$ is a finite set.
\medskip
{\bf Lemma 2}. The set of $z$ for which the homogeneneous polynomial $w \to Hq(w,z)$ vanishes identically is finite.
\medskip
Proof. Write $Hq(w,z)= \sum_{|\alpha|=d} c_\alpha(z)w^\alpha$, and suppose that $H(w,z)$ vanishes identically. After conjugation,
we see that $z$ lies in the variety defined by the vanishing of all the $c_\alpha$.  Consider an index $k$.
Formula (16) shows that each $z_k$ satisfies a polynomial equation with coefficients in $w$ but independent of the other $z_j$.
Hence each $z_k$ can take on at most finitely many values, and the result follows. \ $\spadesuit$
\medskip

{\bf Theorem 4}. Let $q$ be a polynomial that does not vanish on the closed unit ball $B_n$, and let $S(q)$
denote the finite set defined above. If $f={p \over q}$ is a rational proper mapping from $B_n$ to $B_N$ that is reduced to lowest terms,
then $f^{-1}(0) \subset S(q)$. In particular, if $q=1$ then $S(q) = \{0\}$. Thus if
$p$ is a proper polynomial map between balls, then either $p^{-1}(0)$ is empty or $p^{-1}(0)=\{0\}$.

Proof. Given $q$, we suppose that $f={p \over q}$ maps $S^{2n-1}$ to $S^{2N-1}$, and that $f$ is reduced to lowest
terms. By polarization we have

$$ \langle p(z), p(\zeta) \rangle = q(z) {\overline {q(\zeta)}} $$
whenever $\langle z, \zeta \rangle = 1$. Assuming $\langle w,z \rangle \ne 0$, we take $\zeta = {w \over \langle w,z\rangle}$.
We obtain, for all $z$ and $w$ such that $\langle w,z \rangle \ne 0$,

$$ \langle p(z), p({w \over \langle w,z\rangle}) \rangle = q(z) {\overline { q({w \over \langle w,z\rangle})}}. \eqno (17) $$
Assume $||z||<1$ and $p(z)=0$. Then $q(z) \ne 0$ and (17) yield

$$ 0 = \langle w,z \rangle^d q({w \over \langle w,z\rangle})= \sum_{j=0}^{d-1} q_j(w) \langle w,z \rangle^{d-j} = Hq(w,z). \eqno (18) $$
Hence, for each such $z$, the homogeneous polynomial $Hq(w,z)$ vanishes identically.
By Lemma 2, $S(q)$ is a finite set. To prove the second part, we must show that
$S(q) = \{0\}$ when $q$ is the constant polynomial $1$. By construction, in this case $Hq(w,z)= \langle z, w \rangle$, and the only $z$
for which $Hq(w,z)$ vanishes identically in $w$ is $z=0$. Thus either $p^{-1}(0)$ is empty or $p^{-1}(0)=\{0\}$. $\spadesuit$

\medskip
We multiply by $\langle w,z \rangle^d$ in (18) rather than by $\langle w,z \rangle^{d-1}$ in order to take the origin
into account. We next give a different proof of Theorem 4 in the polynomial case; this proof provides additional information.
Homogenization plays the key role in both proofs.

\medskip
For any polynomial map $f$ we let $\nu(f)$ denote its order of vanishing at $0$ and $d(f)$ denote its degree.
Let $p$ be a polynomial proper mapping between balls
and let $p=\sum_{j=0}^d p_j$ denote its expansion into homogeneous parts. Suppose
that $\nu(p) < d(p)$. By reasoning as
in [D] or [D2], the lowest order part $p_\nu$ of $p$ is orthogonal to the highest order part $p_d$;
we may therefore write
$$ p = A \oplus B,$$
where $\nu(A) = \nu(p)$ but $d(A) < d(p)$, and $\nu(B) > \nu(p)$ but $d(B) = d(p)=d$.
The partial tensor 
product operation from [D] replaces $p$ by a proper polynomial map $Ep$, where
$$ Ep = (A \otimes z) \oplus B = A' \oplus B.$$
Evidently $\nu(A')= \nu(A)+1$. It follows that
$\nu(Ep)= \nu(p)+1$ and $d(Ep)=d(p)$. Furthermore, if $p(w)=0$, then $(Ep)(w)=0$ as well and hence the tensor operation
does not decrease the inverse image of the origin.
After iterating the tensor operation a finite number of times, we obtain a homogeneous proper mapping $H$ of degree $d$. By [D]
$H$ must satisfy $||H(z)||^2 = ||z^{\otimes d}||^2= ||z||^{2d}$. Hence the zero set of $H$ is the origin alone. Since
$H(w)=0$ when $p(w)=0$, it follows
that $p^{-1}(0)$ is either empty or just the origin. $\spadesuit$

\medskip

The second proof of Theorem 4 shows that $p$ can be tensored into a map $H$ whose zero set is the origin defined precisely $d$ times.
Furthermore each point in the image of $H$ near $0$ has precisely $d$ inverse images.
For $p$ itself, one cannot assign a single integer mutiplicity to $0$. For example, if $p(z)= (z_1,z_1 z_2,z_2^2)$, then some nearby
points have two inverses images while others have one. Furthermore, for $d\ge 2$,
the juxtaposition of a polynomial of degree $d$ with the identity
provides an injective example of degree $d$. 

Theorem 4 is obvious when the domain and target dimensions both equal $1$; the only proper polynomial mappings are then
already of the form $z^d$, and no tensoring is required. The construction in Proposition 2 is of course analogous
to finding a Blaschke product with a given zero set.

We mention also that Theorem 4 slightly simplifies the proof of the following result in [D1]. If $f$ and $g$ are spherically
equivalent proper polynomial maps which preserve the origin, then they are in fact unitarily equivalent.

\medskip
{\bf V. Valid quadruples}
\medskip

Given a quadruple of integers $(n,r,d,k)$, we naturally ask whether there exists a
$k$-dimensional convex family of proper rational mappings from $B_n$ of generic degree $d$ and rank $r$.
If so we say that $(n,r,d,k)$ is {\it valid}. In this case there is a $k$-dimensional convex family of such maps
with target dimension $N$ whenever $N \ge r$.
Providing necessary and sufficient conditions for a quadruple to be valid
is an extremely difficult problem; in particular, the open problem on degree estimates
is the special case when $k=0$. Nonetheless we give in this section some interesting examples
and constructions of such valid quadruples.

We first make several intuitive statements about valid quadruples.
For a fixed $n$ and $r$ we can make $d$ largest by choosing $k=0$. 
As we increase $k$ we expect that we will decrease $d$. For a fixed $n$ and $d$ the integer $k$ is bounded above; see Proposition 3.
By making $r$ large enough we can choose $k$ to be arbitrarily large as long as $d$ is large enough.
By fixing $r$ we restrict both $d$ and $k$.  These remarks illustrate the basic point. Restrictions on $r$
place restrictions on both measurements of complexity (the degree and the dimension of the family), and furthermore these
measurements are related.

We start with some simple observations about rational proper maps.
For $n\ge 2$, the numerator of a rational example determines the denominator. In other words,
if $||{p \over q_1}||^2 = ||{p \over q_2}||^2$ on the sphere, then $q_1$ is a constant times $q_2$.
The proof is simple; we are given that $|q_1|^2 = |q_2|^2$ on the sphere. Since these functions do not vanish
inside, the maximum principle implies that $q_1$ is a constant times $q_2$. When $n=1$ the conclusion
is false, but an analogous statement holds if we insist that the $q_j$ have no zeroes inside.

The degree of a rational mapping $f$
is the maximum of the degrees of the numerator and denominator when $f$ is reduced to lowest terms.
When $f$ is a proper rational mapping between balls and $f(0)=0$, necessarily its degree is the degree of the numerator. 
See [D] or [F2] for example.
A classical result (in one variable) of Kronecker (see for example [S]) enables one to decide if a given
formal power series defines a rational function. One forms an appropriate Hankel matrix from the Taylor coefficients
and the series defines a rational function if and only if this matrix has finite rank $m$. Furthermore, one can determine
the degree of the rational function from $m$. A similar test works in more variables.

In this section we will provide many examples of valid quadruples $(n,r,d,k)$ and a partial list of restrictions
on them. The families of mappings will be explicit in most cases.

\medskip
{\bf Examples 8}. 

1) $(n,n,1,0)$ is valid for all $n$. Proof: Take $f(z)=z$.

2) Assume $n\ge 2$. Then $(n,r,2,0)$ is invalid if $r < 2n-1$. 
Proof: Faran [Fa1] proved, for $N\le 2n-2$, that a proper rational $f:B_n \to B_N$ has degree at most one.
If $N=2n-1$, then there is a polynomial proper mapping of degree $2$. See [D] or [HJ] for example.
Hence $(n,2n-1,2,0)$ is valid.

3) $(2,a+2,2a+1,0)$ is valid for all $a\ge 0$. Proof: The group invariant monomials used in (11) provide examples. See 
[DKR] and [DLP].

4) $(2,a+4,2a+2,1)$ is valid for all $a\ge 0$. Proof: Take a monomial map $f(z,w)$ of the form
(11). For $\lambda \in (0,1)$ replace the term $|z|^{2(2a+1)}$ in $||f(z,w)||^2$ by 
$$ \lambda |z|^{2(2a+1)} + (1 -\lambda) ( |z|^{2(2a+2)} + |w|^2 |z|^{2(2a+1)}). \eqno (19) $$
The result defines a proper monomial family where the quadruple is
$(2,a+4,2a+2,1)$. We added two terms, increased the degree by one,
and created a one-dimensional convex family. Thus we added $(0,2,1,1)$ to the valid quadruple from 3).
Similar constructions can be used to create higher dimensional families.

\medskip
{\bf Example 9}. Suppose $n=2$ and write $(x_1,x_2) = (x,y)$.
 Consider a monomial example $f$ whose real form $p$ from (12) is of degree $d$
and contains $k$ consecutive monomials of degree $d$.
We write $p = g + \sum_{i=1}^k h_i$, where $h_i$ are these consecutive monomials. The new map defined by
$$ p(x,y) = g(x,y) + \sum_{i=1}^k \lambda_i h_i(x,y) + (x+y) \sum_{i=1}^k(1 -\lambda_i) h_i(x,y)$$
now has $N+k+1$ terms and is of degree $d+1$.  We have created a $k$-dimensional family of monomial maps
at the expense of increasing the degree by $1$ and the target dimension by $k+1$. Thus the quadruple
$(2,N,d,0)$ gets changed into $(2,N+k+1, d+1,k)$. For concreteness we give two explicit
examples. In (20) we have the valid quadruple $(2,5,2,2)$. In (21) we have the valid quadruple $(2,9,3,5)$;
we have written $\lambda = (\lambda_1,...,\lambda_5)$.
$$ \lambda_1 x + \lambda_2 y + (1-\lambda_1)x^2 + (1- \lambda_2)y^2 + (2-\lambda_1 -\lambda_2)xy. \eqno (20)$$
$$ \lambda_1 x + \lambda_2 y + \lambda_3 x^2 + \lambda_4 y^2 + 
\lambda_5 xy + a_1(\lambda) x^3 + a_2(\lambda)x^2 y + a_3(\lambda) xy^2 + a_4(\lambda) y^3 .\eqno (21) $$
Thus (21)
provides a $5$-dimensional parameter space of proper mappings, of generic degree $3$, from $B_2$ to $B_9$.
Because $9$ is the rank of the family, $9$ is also the minimum possible target dimension for this family.

\medskip

We now turn to some general results about valid quadruples.
Let $V(n,d)$ denote the vector space of polynomials of degree at most $d$ in $n$-variables which have no constant term.
Let $\delta(n,d)$ denote its dimension. We next generalize Example 7 by finding the largest possible dimension
for a family of origin preserving polynomial maps of degree $d$.

\medskip
{\bf Proposition 3}. Fix $n$ and $d$. Let $k = \delta(n,d-1)^2$.
Then there is a $k$-dimensional convex family
of proper polynomial mappings from $B_n$ of generic degree $d$ and of rank $\delta(d,n)$.
Furthermore, for every convex family of origin preserving polynomial proper mappings from $B_n$ of degree $d$, 
the dimension of the family is at most this number $k$.

Proof. We begin by considering a Hermitian symmetric polynomial $R$ of degree $d$ and vanishing at $0$.
Thus 
$$ R(z, {\overline z}) = \sum_{1 \le |\alpha|, |\beta| \le d} c_{\alpha \beta} z^\alpha {\overline z}^\beta, \eqno (22) $$
where we regard the coefficients as parameters. Let $C$ denote the Hermitian matrix of coefficients.
We consider the system of equations for these parameters obtained by setting $R(z,{\overline z})=1$ on the sphere.
Reasoning as in Example 7 we discover that these parameters satisfy a universal system of linear equations; we can
always solve for those $c_{\alpha \beta}$ where at least one of $|\alpha|$ or $|\beta|$ equals $d$, in terms of those
where both multi-indices are of order less than $d$. We let $A$ denote the Hermitian submatrix of $C$ defined by
those $c_{\alpha \beta}$ where $|\alpha|$ and $|\beta|$ are at most $d-1$. We can therefore write 
$$ C = \left ( \matrix {  A & B
 \cr  B^*  &  D } \right), $$
where $B$ and $D$ are determined by $A$. In fact the entries of $B$ are linear functions of the entries in $A$,
whereas those in $D$ are affine functions of the entries in $A$. When $A=0$, we must have $B=0$, and
these equations have a solution where $D_0(z,{\overline z}) = ||z||^{2d}$. The eigenvalues of the corresponding diagonal
form $D_0$ are the multinomial coefficients, and hence $D_0$ is positive definite on its domain. We may choose the matrix $A$ as we wish,
and hence we assume that it is positive definite on its domain. Since positive definiteness is an open condition,
it follows by the continuity of linear and affine functions that
the full matrix $C$ will be positive definite if $A$ is positive definite and its entries are sufficiently small. The entries
of $B$ are forced to be small when the entries of $A$ are. Thus there is an
open subset in the parameter space for which $C$ is positive definite.
In this open set it follows that $R(z,{\overline z}) = ||f(z)||^2$ for some polynomial $f$ of degree $d$ depending 
on the parameters from $A$; as these parameters vary we obtain
a $k$-dimensional family, where $k$ is the dimension over ${\bf R}$ of the space of Hermitian forms on $V(n,d-1)$.
This dimension is $\delta(n,d-1)^2$. The rank of a family is the rank of its underlying Hermitian form; since the matrix $C$
is positive definite, it is invertible, and hence its rank is $\delta(n,d)$. 

In our construction we allowed the maximum possible number of parameters for the squared norm
of an origin preserving proper polynomial map of degree $d$ in $n$ variables. We showed that we can always construct
a family whose dimension realizes this maximum; the last statement follows.
\ $\spadesuit$
\medskip

In the proof of Proposition 3 we assumed that the map $f$ was of degree $d$, and we showed that its Taylor polynomial
of degree $d-1$ at the origin essentially determined $f$. Given the Taylor polynomial of degree $d-1$, bounds on the target dimension
of the type we have considered might force $f$ to be of degree $d$. In such cases $f$ is essentially determined by its $d-1$
jet at the origin. See [EL] for finite jet determination results for CR imbeddings of Levi nondegenerate hypersurfaces. In [EL]
the map must satisfy a quantitative nondegeneracy condition at the base point. The jet needed to determine the map depends on this
quantity. It is natural to wonder how our results on valid quadruples relate to finite jet detrmination in a more general setting.

\medskip
We next use Example 7 and Proposition 3 when $d=2$ to indicate how complicated the set $K$ from Lemma 1 can be. 
Consider four real parameters given by real variables $x,y$ and a complex variable $\zeta$. We start with a Hermitian form
$A$ in two variables:
$$ A = \left ( \matrix {  x & \zeta
 \cr  {\overline \zeta}  &  y} \right). $$
The matrix $B$ from the proof of Proposition 3 will be the zero matrix. The proof of Proposition 3 shows that
the matrix $D$ is given by
$$ D= \left ( \matrix {  1-x & -\zeta & 0
 \cr  -{\overline \zeta} & 2-x-y  &  -\zeta \cr 0 & -{\overline \zeta} & 1-y } \right). $$
Because $B=0$, the condition for $C$ to be a nonnegative definite (five by five) matrix is that
both $A$ and $D$ are nonnegative definite. The three conditions for $C$ being nonnegative definite
turns out to be $0\le x,y \le 1$ and $|\zeta|^2 \le xy$ and $|\zeta|^2 \le (1-x)(1-y)$. The boundary of this set
is evidently somewhat complicated.

\medskip 
We mention an alternative viewpoint for Proposition 3.
Let ${\cal H}$ be a complex Hilbert space with inner product $\langle, \rangle$. Consider
a polynomial function $p: {\bf C}^n \to {\cal H}$ of degree $d$ such that $p(0)=0$. Thus
$$ p(z) = \sum_{|\alpha|=1}^d C_\alpha z^\alpha, \eqno (23) $$
where each $C_\alpha \in {\cal H}$.
The condition that $||p(z)||^2 = 1$ on $||z||^2=1$ is a system of finitely many linear equations
in the finitely many unknowns given by the inner products $\langle C_\alpha, C_\beta \rangle$. These are the parameters $c_{\alpha \beta}$
from Proposition 3. As before we can solve this system of equations for the inner products
$\langle C_\alpha, C_\beta \rangle$ where one or both
of $|\alpha|$ or $|\beta|$ equals $d$, in terms of those inner products indexed by the multi-indices
of order less than $d$. It follows that the $d-1$ jet of $p$ determines $p$, up to norm equivalence.
Thus the dimension of the parameter space equals the dimension over ${\bf R}$ of the space of Hermitian forms on $V(n,d-1)$.

The next result appears in [D]; we state it here and sketch its proof
to clarify the subtle contrast between it and Proposition 3. In Proposition 3 the parameters from the matrix $C$ need
to be small; the condition differs from the condition that the function $R$ map the closed ball to the open ball. In this latter case,
we can control neither the degree nor the rank (embedding dimension).

\medskip
{\bf Proposition 4}. Let $g:{\bf C}^n \to {\bf C}^N$ be a polynomial mapping with $g(0)=0$ and $||g(z)|| < 1$ on the closed ball.
Then there is an integer $K$ and polynomial mapping $p:{\bf C}^n \to {\bf C}^K$ such that $p(0)=0$
and $g \oplus p$ maps $B_n$ properly to $B_{N+K}$. No bound on either $K$ or the degree of $p$ depending
on only $n$ and the degree of $g$ exists.

Proof. (Sketch) We suppose that $g$ is of degree $d-1$ to maintain the parallel with Proposition 3.
Reasoning as in [D2] we proceed as follows. Let $R(z,{\overline z}) = ||z||^{2(d-1)} - ||g(z)||^2$.
Then $R$ defines a Hermitian symmetric polynomial which is strictly positive on the unit sphere. By
[CD1] there is a polynomial mapping $p:{\bf C}^n \to {\bf C}^K$ such that
$$ R(z,{\overline z}) = ||p(z)||^2 \eqno (24) $$
on the unit sphere and we may assume that $p(0)=0$. Then $g \oplus p$ does the job.

The existence of $p$ satisfying (24) is non-trivial, and no bounds on its degree or on $K$ in terms
of $n$ and $d-1$ are possible. See [CD1] or [D3].
One first adds a variable to bihomogenize $R$, and one adds a term
to guarantee that the bihomogenized $R$, written $R_h$,  is strictly positive on the unit sphere in ${\bf C}^{n+1}$. 
This condition does not imply the result. One must multiply by a sufficiently high power of $||z||^2 + |t|^2$ to guarantee 
that the underlying form of
$$ (||z||^2 + |t|^2)^m R_h(z,t,{\overline z}, {\overline t}) \eqno (25) $$
is positive definite, and hence of the form $||P(z,t)||^2$. One then obtains (24) by first dehomogenizing (25) and 
then restricting to the unit sphere. \ $\spadesuit$

\medskip

We close this paper with two remarks.
The first concerns the distinction between the monomial case and the polynomial case.
Consider the proof of Proposition 3. When $f$ is a monomial, all off-diagonal elements of the underlying 
matrix vanish. We are therefore imposing additional constraints that lower the possible dimension of a family.

The sketch of the proof of Proposition 4 illustrates a crucial point. Let $R(z, {\overline z})$ be a bihomogeneous Hermitian
symmetric polynomial. The condition for $R$ to be a squared norm is stronger than its non-negativity. When
$R$ is nonnegative, we cannot even conclude that it is a quotient of squared norms. If, however, $R$ 
is strictly positive on the unit sphere, then there is an integer $m$ such that 
$||z||^{2m} R(z, {\overline z})$  is a squared norm. 
It follows that $R$ agrees with a squared norm on the sphere.
See [D2],[D3],[CD1], [CD2] for lengthy discussions of this point and its generalization to isometric imbedding theorems.

\medskip

\medskip
{\bf References}
\medskip

[CD1] Catlin, D. and D'Angelo, J., A stabilization theorem for 
Hermitian forms and applications to holomorphic mappings, {\it Math. Res. Lett.} 3 (1996), 149-166.
\medskip

[CD2] Catlin, D. and D'Angelo, J., An isometric imbedding theorem for holomorphic bundles,
{\it Math. Res. Lett.} 6 (1999), 43-50.
\medskip

[D] D'Angelo, J., Several Complex Variables and the Geometry of Real Hypersurfaces,
CRC Press, Boca Raton, 1993.
\medskip

[D1] D'Angelo, J., Proper polynomial mappings between balls of different dimensions, {\it Michigan Math J.}
35 (1988), 83-90.

\medskip
[D2] D'Angelo, J., Proper holomorphic mappings, 
positivity conditions, and isometric imbedding, {\it J. Korean Math Society}, May 2003, 1-30.

\medskip
[D3] D'Angelo, J., Inequalities from Complex Analysis, Carus Mathematical Monograph No. 28, 
Mathematics Association of America, 2002.

\medskip
[D4] D'Angelo, J., Homogenization, Reflection, and the $X$-variety, {\it Indiana Univ. Math J.} 52 (2003),
1113-1134.

\medskip

[DKR] D'Angelo, J., Kos, \v S., and Riehl, E.,
A Sharp Bound for the Degree of Proper Monomial
Mappings Between Balls, {\it J. Geometric Analysis}, Volume 13 (2003), no. 4,
581-593. 

\medskip 

[DLP] D'Angelo, J., Lebl, J.,  and Peters, H.,
Degree estimates for polynomials constant on hyperplanes,  {\it Michigan Math. J.}, 55 (2007), no. 3, 693-713.
\medskip

[DV] D'Angelo, J. and Varolin, D.,  Positivity conditions for
Hermitian symmetric functions, {\it Asian J. Math} 8 (2004), 215-232.

\medskip
[Dor] Dor, A., Proper holomorphic maps between balls in one co-dimension, {\it  Ark. Mat.}  28  (1990),  no. 1, 49-100. 
\medskip
[EHZ] Ebenfelt, P., Huang, X.,  and Zaitsev, D., Rigidity of CR-immersions into spheres,
{\it Comm. Anal. Geom.},  12  (2004),  no. 3, 631-670.

\medskip
[EL] Ebenfelt, P. and Lamel, B., Finite Jet Determination of CR Embeddings,
 {\it J. Geometric Analysis}, Volume 14 (2004), no. 2, 241-265.

\medskip
[Fa1] Faran, J., On the linearity of 
proper maps between balls in the low codimension case,
{\it J. Diff. Geom.}, 24 (1986), 15-17.
\medskip
[Fa2] Faran, J., Maps from the two-ball to the three-ball, 
{\it Inventiones Math.}, 68 (1982), 441-475.

\medskip
[F1] Forstneric, F., 
Extending proper holomorphic maps of positive codimension, 
{\it Inventiones Math.}, 95(1989), 31-62.
\medskip

[F2] Forstneric, F., Proper rational maps: A survey, Pp 297-363 in
{\it Several Complex Variables: Proceedings of the Mittag-Leffler
Institute, 1987-1988},  Mathematical Notes 38, Princeton Univ.
Press, Princeton, 1993.

\medskip
[H] Huang, X., On a linearity problem for proper maps between balls in complex spaces
of different dimensions, {\it J. Diff. Geometry} 51 (1999), no 1, 13-33.

\medskip
[HJ] Huang, X., and Ji, S., 
Mapping $B_n$ into $B_{2n-1}$, {\it Invent. Math.} 145 (2001), 219-250.
13-36.
\medskip 
[HJX] Huang, X.,  Ji, S.,  and Xu,D.,
Several results for holomorphic mappings from ${\bf B}_n$ to 
${\bf B}_N$, {\it Contemporary Math.} 368 (2005), 267-292.

\medskip

[HJX2] Huang, X., Ji, S.,  and D. Xu, 
A new gap phenomenon for proper holomorphic mappings from $B\sp n$ into $B\sp N$.
{\it Math. Res. Lett.} 13 (2006), no. 4, 515--529. 

\medskip
[L] Lebl, J.,
Singularities and Complexity in CR Geometry, PhD thesis, 
University of California, San Diego, 2007.

\medskip

[M] Meylan, F., Degree of a holomorphic map between unit balls from ${\bf C}^2$ 
to ${\bf C}^n$,  {\it Proc. Amer. Math. Soc.}  134  (2006),  no. 4, 1023-1030.

\medskip
[P] Pinchuk, S. Holomorphic maps in ${\bf C}^n$ and the problem of holomorphic
equivalence. Pp 173-200 in {\it Encyclopedia of Math. Sciences}, 
Vol. 9, G. M. Henkin, editor,
Several Complex Variables III, 173, 1980.

\medskip
[S] Sarason, D., review of Hankel operators and their applications by Vladimir V. Peller,
{\it Bull. A. M. S.}, Volume 41, Number 3, (2004) 401-407.

\end